\documentclass[11pt,leqno,fleqn]{article}
\usepackage{amsmath,amssymb}
\newcounter{sectie}
\newcommand{\sectz}[1]{\refstepcounter{sectie}
\section*{\boldmath \thesectie. #1}%
}
\newcommand{\dy}[2]{%
\refstepcounter{equation}%
\label{#1}%
\begin{list}{}{
\topsep 3mm
\leftmargin 18mm
\rightmargin 0cm
\itemsep 0mm
\listparindent 0mm
\parsep 0mm
\itemsep 0mm
\labelsep 0mm
\labelwidth 18mm
}%
\item[\rm (\theequation)\hfill]
#2
\end{list}%
}
\newcommand{\diz}[1]{%
\refstepcounter{equation}%
\begin{list}{}{
\topsep 5mm
\leftmargin 10mm
\rightmargin 0cm
\itemsep 0mm
\listparindent 0mm
\parsep 0mm
\labelsep 1mm
\labelwidth 10mm
}%
\item[\rm (\theequation)\hfill]
\begin{list}{}{
\topsep 0mm
\leftmargin 8mm
\rightmargin 0mm
\itemsep 0mm
\listparindent 0mm
\parsep 0mm
\labelsep 1.5mm
\labelwidth 6.5mm
}
#1
\end{list}%
\end{list}%
}
\newcommand{\nr}[1]{\item[{\rm (#1)}]}
\newcommand{\nrs}[1]{\item[{\rm (#1)}]\vspace{-\itemsep}}
\newcommand{\oR}{{\mathbb{R}}}
\newcommand{\corank}{{\text{\rm corank}}}
\newcommand{\dez}[1]{\dyz{\raggedright$\displaystyle #1 $}}
\newcommand{\dyz}[1]{%
\refstepcounter{equation}%
\begin{list}{}{
\topsep 3mm
\leftmargin 18mm
\rightmargin 0cm
\itemsep 0mm
\listparindent 0mm
\parsep 0mm
\itemsep 0mm
\labelsep 0mm
\labelwidth 18mm
}%
\item[\rm (\theequation)\hfill]
#1
\end{list}%
}
\newcounter{bewering}
\newcommand{\propz}[1]{\refstepcounter{bewering}\vspace{4mm}\noindent{\bf Proposition \thebewering.}{\it #1}}
\newcommand{\pf}{\vspace{3mm}\noindent{\bf Proof.}\ }
\newcommand{\T}{^{\sf T}}
\newcommand{\rf}[1]{{\rm (\ref{#1})}}
\newcommand{\bx}{\hspace*{\fill} \hbox{\hskip 1pt \vrule width 4pt height 8pt depth 1.5pt \hskip 1pt}

\addvspace{4mm}}
\newcommand{\prop}[2]{\refstepcounter{bewering}\vspace{4mm}\noindent{\bf Proposition \thebewering.}\label{#1}{\it #2}}
\newcommand{\propnmz}[2]{\refstepcounter{bewering}\vspace{4mm}\noindent{\bf Proposition \thebewering\ {\rm (#1)}.}{\it #2}}
\newcommand{\di}[2]{%
\refstepcounter{equation}%
\label{#1}%
\begin{list}{}{
\topsep 5mm
\leftmargin 10mm
\rightmargin 0cm
\itemsep 0mm
\listparindent 0mm
\parsep 0mm
\labelsep 1mm
\labelwidth 10mm
}%
\item[\rm (\theequation)\hfill]
\begin{list}{}{
\topsep 0mm
\leftmargin 8mm
\rightmargin 0mm
\itemsep 0mm
\listparindent 0mm
\parsep 0mm
\labelsep 1.5mm
\labelwidth 6.5mm
}
#2
\end{list}%
\end{list}%
}
\newcommand{\PP}{{\cal P}}
\newcommand{\LL}{{\cal L}}
\newcounter{hulpstelling}
\newcommand{\lemma}[2]{\refstepcounter{hulpstelling}\vspace{4mm}\noindent{\bf Lemma \thehulpstelling.}\label{#1}{\it #2}}
\newcommand{\dyyz}[1]{\dyz{\raggedright$\dps#1$}}
\newcommand{\dps}{\displaystyle}
\newcommand{\thmnn}[1]{\vspace{4mm}\noindent{\bf Theorem.}{\it #1}}
\begin{document}

\begin{center}
{\large\bf THE STRONG ARNOLD PROPERTY FOR 4-CONNECTED FLAT GRAPHS

}
\vspace{4mm}
Alexander Schrijver\footnote{ Korteweg-de Vries Institute for Mathematics, University of Amsterdam}
and
Bart Sevenster$\mbox{}^1$

\end{center}

\noindent
{\small
{\bf Abstract.}
We show that if $G=(V,E)$ is a 4-connected flat graph, then any
real symmetric $V\times V$ matrix $M$ with exactly one negative eigenvalue and
satisfying, for any two distinct vertices $i$ and $j$,
$M_{ij}<0$ if $i$ and $j$ are adjacent, and
$M_{ij}=0$ if $i$ and $j$ are nonadjacent,
has the Strong Arnold Property:
there is no nonzero real symmetric $V\times V$ matrix $X$
with $MX=0$ and $X_{ij}=0$ whenever $i$ and $j$ are equal or adjacent.
(A graph $G$ is {\em flat} if it can be embedded injectively
in $3$-dimensional Euclidean space such that the image of any circuit is the boundary of some disk disjoint
from the image of the remainder of the graph.)

This applies to the Colin de Verdi\`ere graph parameter, and extends similar results
for 2-connected outerplanar graphs and 3-connected planar graphs.

}



\sectz{Introduction}

Let $G=(V,E)$ be an undirected graph.
Call a real symmetric $V\times V$ matrix a {\em well-signed $G$-matrix} if for all
distinct $i,j\in V$: $M_{ij}<0$ if $ij\in E$ and $M_{ij}=0$ if $ij\not\in E$.
(No condition on the diagonal elements.)
For any real symmetric matrix $M$, let $\lambda^-(M)$ be the number of negative
eigenvalues, taking multiplicities into account.
The {\em corank} of $M$ is the dimension of its nullspace $\ker(M)$.

The famous {\em Colin de Verdi\`ere parameter} $\mu(G)$ [1] is defined to be the maximal corank of any well-signed $G$-matrix $M$
with $\lambda^-(M)=1$ and having the {\em Strong Arnold Property}:
\dy{21no15b}{
there is no nonzero real symmetric $V\times V$ matrix $X$ with $MX=0$ and $X_{ij}=0$ whenever
$i$ and $j$ are equal or adjacent.
}

The interest of the parameter $\mu(G)$ was exhibited by Colin de Verdi\`ere [1], who showed that
$\mu(G)$ is {\em minor-monotone}, that is, $\mu(H)\leq\mu(G)$ if $H$ is a minor of $G$,
--- in other words, for each $k$, the collection of graphs $G$ with $\mu(G)\leq k$ is closed under
taking minors; hence there are finitely many forbidden minors, by Robertson and Seymour [9].
The Strong Arnold Property is crucial for the minor-monotonicity.

Moreover, Colin de Verdi\`ere [1] showed (i), (ii), and (iii) in:
\diz{
\nr{i} $\mu(G)\leq 1$ if and only if $G$ is a disjoint union of paths,
\nrs{ii} $\mu(G)\leq 2$ if and only if $G$ is outerplanar,
\nrs{iii} $\mu(G)\leq 3$ if and only if $G$ is planar,
\nrs{iv} $\mu(G)\leq 4$ if and only if $G$ is flat.
}
Statement (iv) was proved by Robertson, Seymour, and Thomas [10] (only if) and Lov\'asz and Schrijver [6] (if).
Recall that a graph $G$ is {\em flat} if it can be embedded injectively
in $\oR^3$ such that the image of any circuit is the boundary of some disk disjoint
from the image of the remainder of the graph.
As was shown in [10],
a graph is flat if and only if it is linklessly embeddable, that is, can be
embedded injectively in $\oR^3$ such that the images of any two disjoint circuits
are unlinked.
We refer to [4]
for a survey of the Colin de Verdi\`ere parameter.

A short proof of (iii) was given by van der Holst [2], which proof also implies
that if $G$ is 3-connected and planar, then any well-signed $G$-matrix $M$ with
$\lambda^-(M)=1$, has corank at most 3.
So the Strong Arnold Property is not needed to define $\mu(G)$ for such graphs $G$.
That is, if we define $\kappa(G)$ to be the maximum corank of any well-signed
$G$-matrix $M$ with $\lambda^-(M)=1$,
then $\kappa(G)=\mu(G)$ for 3-connected planar graphs $G$.
Here 3-connectivity cannot be relaxed to 2-connectivity, since
$\kappa(K_{2,t})=t$ for all $t$, while $\mu(K_{2,t})=3$ for all $t\geq 3$.
In [6], it was shown that $\kappa(G)=\mu(G)$ also for 
4-connected flat graphs.

The latter means that for any 4-connected flat graph $G$,
among the well-signed $G$-matrices $M$ with $\lambda^-(M)=1$ that
maximize $\corank(M)$,
there is one having the Strong Arnold Property.
In this paper, we prove that for any 4-connected flat graph $G$,
{\em each} well-signed $G$-matrix $M$ with $\lambda^-(M)=1$
has the Strong Arnold Property.
This extends results of van der Holst [3] who proved this for 2-connected outerplanar
graphs and for 3-connected planar graphs.
In fact, one may show that if this holds for all $\mu(G)$-connected graphs with
$\mu(G)=k$, then also for all $\mu(G)$-connected graphs $G$ with $\mu(G)\leq k$
(by an apex graph argument).

The above raises the question whether the Strong Arnold Property would be superfluous
to impose for all $\mu(G)$-connected graphs $G$ --- in the weak sense: that
$\kappa(G)=\mu(G)$, or in the strong sense: that each well-signed $G$-matrix $M$
with $\lambda^-(M)=1$, has the Strong Arnold Property.
We do not put this as conjecture, since our proof method might suggest that the
case $\mu(G)\leq 4$ is exceptional.

The relevance of the present paper may also lie in obtaining a better understanding
of the nullspace embedding of a graph $G$ defined by $M$ (see below).
For a 3-connected planar graph, such a nullspace embedding corresponds to a planar embedding
of the graph on the 2-sphere (Lov\'asz [5], cf.\ [7,\linebreak[0]8]).
An intriguing question is whether, if $G$ is a 4-connected flat
graph and $M$ is a well-signed $G$-matrix with $\lambda^-(M)=1$,
its nullspace embedding (normalized to unit-length vectors) yields a flat embedding of $G$ on the
3-sphere.
The fact that any such matrix has the Strong Arnold Property may help in proving this.

\sectz{The Strong Arnold Property and quadrics}

We first formulate the Strong Arnold Property of $M$ in terms of the nullspace embedding
defined by $M$.
Let $G=(V,E)$ be an undirected graph and let $M$ be a well-signed $G$-matrix
with $\lambda^-(M)=1$ and with corank $d$.
Let $b_1,\ldots,b_d\in\oR^V$ be a basis of $\ker(M)$.
Define, for each $i\in V$, the vector $u_i\in\oR^d$ by: $(u_i)_j:=(b_j)_i$, for
$j=1,\ldots,d$.
So we have $u:V\to\oR^d$.
Then $u$ is called the {\em nullspace embedding of $G$ defined by $M$}.
Note that $u$ is unique up to linear transformations of $\oR^d$.

The Strong Arnold Property of $M$ is in fact a property only of the graph
$G$ and the function $i\mapsto\langle u_i\rangle$.
(Throughout, $\langle\ldots\rangle$ denotes the linear space spanned by $\ldots$.)
When we have $u:V\to\oR^d$, define $|G|$ to be the following subset of $\oR^d$:
\dez{
|G|:=\bigcup\{\langle u_i\rangle\mid i\in V\}
\cup
\bigcup\{\langle u_i,u_j\rangle\mid ij\in E\}.
}
A subset $Q$ of $\oR^d$ is called a {\em homogeneous quadric} if it is the solution set
of a nonzero homogeneous quadratic equation.

\propz{
$M$ has the Strong Arnold Property if and only $|G|$ is not contained in any homogeneous
quadric.
}

\pf
Let $U$ be the $d\times V$ matrix with as columns the vectors $u_i$ for $i\in V$.

Suppose that some homogeneous quadric $Q=\{y\mid y\T Ny=0\}$ contains $|G|$,
where $N$ is a nonzero symmetric $d\times d$ matrix.
Then $X:=U\T NU$ is a nonzero symmetric $V\times V$ matrix that contradicts
the Strong Arnold Property \rf{21no15b}.

Conversely, suppose that $M$ has not the Strong Arnold Property.
Let $X$ be a matrix as in \rf{21no15b}.
As $MX=0$ and as $X$ is symmetric, we have $X=U\T NU$ for some nonzero symmetric
$d\times d$ matrix $N$.
Then $Q:=\{y\mid y\T Ny=0\}$ is a homogeneous quadric containing $|G|$.
\bx

Throughout this paper,
by a hyperplane, plane, and line in $\oR^d$ we mean {\em linear} subspaces, of dimension
$d-1$, $2$, and $1$, respectively.
Note that if a homogeneous quadric $Q$ contains a hyperplane then
$Q$ is the union of one or two hyperplanes
(as we can assume that $Q$ contains $\{x\mid x_1=0\}$,
hence the quadratic form is $(a\T x)x_1$ for some nonzero $a\in\oR^d$).

We will consider triples $G,M,u$ where
\dy{26no15a}{
$G$ is a graph,
$M$ is a well-signed $G$-matrix with one negative eigenvalue,
and $u:V\to\oR^d$ is the nullspace embedding defined by $M$.
}
The essence of our proof is showing that, for any such triple $G,M,u$
with $G$ a 4-connected flat graph and $d=4$,
\dy{27no15b}{
$|G|$ is not contained in the union of two hyperplanes, and
$|G|$ contains distinct planes $P_1,\ldots,P_4$ with $P_1\cap P_2\cap P_3\neq\{0\}$ and
$P_1\cap P_4=\{0\}$.
}
Having this, the following basic fact on quadrics shows that $|G|$ cannot be contained in
any homogeneous quadric:

\prop{25no15a}{
Let $Q$ be a homogeneous quadric in $\oR^4$, not being the union of two hyperplanes.
If a line is contained in three planes on $Q$, it is contained in each plane on $Q$.
}

\pf
Suppose line $\ell$ is contained in planes $P_1,P_2,P_3$ on $Q$, but not in plane
$R$ on $Q$.
Consider two distinct $i,j\in\{1,2,3\}$, and define $H:=P_i+P_j$.
As $H\not\subseteq Q$, $Q':=Q\cap H$ is a homogeneous quadric in $H$.
Since $Q'\supseteq P_i\cup P_j$, we know $Q'= P_i\cup P_j$.
Hence, as $R\cap H\neq\{0\}$ (since $\dim(R)=2$ and $\dim(H)=3$)
and as $P_i\cap P_j\cap R=\ell\cap R=\{0\}$, $R\setminus\{0\}$
intersects precisely one of $P_i$ and $P_j$.
As this cannot hold simultaneously for each two $i,j$ in $\{1,2,3\}$, we are done.
\bx

\sectz{Graphs and hyperplanes}

Having $u:V\to\oR^d$, we say that a subspace is {\em spanned} if it is linearly
spanned by a subset of $u(V)$.
A crucial tool will be the following lemma of van der Holst [2]:

\propnmz{Van der Holst's lemma}{
Let $G,M,u$ satisfy \rf{26no15a}, and let
$H$ be a hyperplane in $\oR^d$, splitting $\oR^d$ into the two halfspaces
$H'$ and $H''$.
\di{27no15d}{
\nr{i} If $G$ is connected and $H$ is spanned, then 
each of the vertex sets $u^{-1}(H')$ and $u^{-1}(H'')$ is nonempty
and spans a connected subgraph of $G$;
\nrs{ii}
Any vertex in $u^{-1}(H)$ with a neighbour in $u^{-1}(H')$, has also a 
neighbour in $u^{-1}(H'')$.
}
}

Van der Holst's lemma gives the first half in \rf{27no15b}:

\prop{18no15a}{
Let $G,M,u$ satisfy \rf{26no15a}, with $G$ a $4$-connected flat graph.
Then $|G|$ is not contained in the union of two hyperplanes.
}

\pf
Suppose $|G|\subseteq H_1\cup H_2$ for hyperplanes $H_1$, $H_2$ in $\oR^d$.
As $u(V)$ is full-dimensional, $H_1$ and $H_2$ are distinct, and we can assume
they are spanned hyperplanes.
Let $H'_i$ and $H''_i$ be the two sides of $H_i$.
By \rf{27no15d}(i), for each $i=1,2$,
each of the vertex sets $u^{-1}(H'_i)$ and $u^{-1}(H''_i)$ induces a connected
subgraph of $G$.
As $G$ is $4$-connected, there exist $4$ internally disjoint paths connecting
$u^{-1}(H_1')$ and $u^{-1}(H_2')$.
By \rf{27no15d}(ii), $H_1'$ and $H_1''$ have the same neighbours
in $u^{-1}(H_1\cap H_2)$.
Similarly, $H_2'$ and $H_2''$ have the same neighbours in $u^{-1}(H_1\cap H_2)$.
Hence we can assume that all internal vertices of these paths belong to
$u^{-1}(H_1\cap H_2)$.
Contracting each of these paths, and contracting each $u^{-1}(H'_i)$ and
$u^{-1}(H''_i)$, we obtain $K_{4,4}$.
This is a contradiction, as $K_{4,4}$ is not flat.
\bx

From Proposition \ref{18no15a} we derive:

\prop{26no15bx}{
Let $G,M,u$ satisfy \rf{26no15a}, with $G$ a $4$-connected flat graph and $d=4$.
Then there exist planes $P,R\subseteq|G|$ with $P\cap R=\{0\}$.
}

\pf
Let $\PP$ be the collection of planes $P\subseteq|G|$, and let $\LL$ be the
collection of lines in $|G|$ not contained in any plane.
Suppose to the contrary that $P\cap R\neq\{0\}$ for all $P,R\in\PP$.

Let $H$ be a spanned hyperplane containing a maximum number of planes in $\PP$.
If some $P\in\PP$ is not contained in $H$, then there exist $R,S\in\PP$ with
$H=R+S$ (by the maximality), and $R\cap S\subseteq P$
(as $R\cap P\neq\{0\}$ and $S\cap P\neq \{0\}$, while $P\not\subseteq R+S$).
Hence there exists a line $\ell\subset H$ with $P\cap H=\ell$ for each
$P\in \PP$ with $P\not\subseteq H$.
Concluding, for all distinct $P,R\in\PP\cup\LL$, $P\setminus H$ and $R\setminus H$ are
disjoint (as if $P,R\in\PP$ then $\ell=P\cap R$, hence $P\setminus\ell$ and
$R\setminus\ell$ are disjoint).

Let $H'$ and $H''$ be the halfspaces separated by $H$.
By \rf{27no15d}(i), $u^{-1}(H')$ and $u^{-1}(H'')$ induce connected subgraphs
 of $G$.
So there are at most two $P\in\PP\cup\LL$ with $P\not\subseteq H$.
Let $J$ be the sum of these $P$.
Then $J$ has dimension at most 3.
So $|G|$ is contained in the union of two hyperplanes, contradicting
Proposition \ref{18no15a}.
\bx

\sectz{Existence of $P_1,P_2,P_3,P_4$}

In this section, we prove the second half in \rf{27no15b}.
First, three lemmas.

\lemma{18no15j}{
Let $G,M,u$ satisfy \rf{26no15a}, with $d=\kappa(G)\geq 2$ and $G$ connected.
Let $G'$ be a subgraph of $G$ with $V(G')=V$, and let $A$ be a well-signed $G'$-matrix
with $\ker(M)\subseteq\ker(A)$.
Then $\lambda^-(A)\leq 1$.
}

\pf
Suppose $\lambda^-(A)\geq 2$.
Then there exists $\beta>0$ such that $\lambda^-(\beta A+M)\geq 2$.
Let $\alpha$ be the infimum of these $\beta$.
Note that $\corank(\beta A+M)\geq\corank(M)$, since $\ker(M)\subseteq\ker(A)$.
For any real symmetric matrix $X$, denote by $\lambda_i(X)$ the $i$-th eigenvalue of $X$ from
below, taking multiplicities into account.

Then $\lambda^-(\alpha A+M)=1$.
Suppose not.
Then $\alpha>0$.
As $\alpha A+M$ is a $G$-matrix, as $G$ is connected, and as
$\corank(\alpha A+M)\geq\corank(M)\geq 2$,
$\alpha A+M$ has at least one negative eigenvalue, by Perron-Frobenius.
For each $\gamma<\alpha$ one has $\lambda^-(\gamma A+M)\leq 1$, so
$\lambda_2(\gamma A+M)\geq 0$, hence by continuity of $\lambda_2$,
$\lambda_2(\alpha A+M)\geq 0$.
So $\lambda^-(\alpha A+M)=1$, contradicting our assumption.

Moreover, by definition of $\alpha$, there exist $\beta>\alpha$ arbitrarily close to $\alpha$ with $\lambda^-(\beta A+M)\geq 2$.
As $\corank(\beta A+M)\geq\corank(M)=\kappa(G)=:k$, we have $\lambda_{k+2}(\beta A+M)\leq 0$.
Then, by continuity of $\lambda_{k+2}$, $\lambda_{k+2}(\alpha A+M)\leq 0$.
So $\corank(\alpha A+M)\geq k+1>\kappa(G)$.
This contradicts the definition of $\kappa(G)$.
\bx

\lemma{18no15i}{
Let $C$ be a circuit and let $u:V(C)\to\oR^2\setminus\{0\}$ be such that
for any two incident edges $ij$ and $jk$,
the vectors $u_i$ and $u_k$ are at different sides of the line $\langle u_j\rangle$.
Then there exists a well-signed $C$-matrix $A$ with $\lambda^-(A)\geq 1$ such that
$u$ is the nullspace embedding defined by $A$.
}

\pf
For each edge $ij$ of $C$, define $a_{ij}:=-|\det(u_i,u_j)|^{-1}$.
If $i$ and $k$ are the two neighbours of vertex $j$,
then $v:=a_{ij}u_i+a_{jk}u_k$ is a scalar multiple of $u_j$;
equivalently, $\det(v,u_j)=0$.
Indeed,
\dyyz{
\det(v,u_j)=
a_{ij}\det(u_i,u_j)+a_{jk}\det(u_k,u_j)
=
-\frac{\det(u_i,u_j)}{|\det(u_i,u_j)|}
-
\frac{\det(u_k,u_j)}{|\det(u_j,u_k)|}
=0,
}
since $\det(u_i,u_j)$ and $\det(u_k,u_j)$ have opposite signs
(as $u_i$ and $u_k$ are at different sides of $\langle u_j\rangle$).

Concluding, there exists $a_{jj}$ such that $v=-a_{jj}u_j$, yielding the matrix $A$.
Note that necessarily $\lambda^-(A)\geq 1$, by Perron-Frobenius, as $\corank(A)\geq 2$
and $C$ is connected.
\bx

\lemma{18no15h}{
Let $G,M,u$ satisfy \rf{26no15a}.
Let $P\subseteq|G|$ be a plane such that there are no two other planes
$R,S\subseteq|G|$ with $P\cap R\cap S\neq\{0\}$.
Then there exists a subgraph $G_P$ of $G$ which is a circuit on a subset
of $u^{-1}(P\setminus\{0\})$ added with
isolated vertices, and a well-signed $G_P$-matrix $A_P$
with $\ker(M)\subseteq\ker(A_P)$ and $\lambda^-(A_P)\geq 1$.
}

\pf
Choose an edge $ij$ such that $P=\langle u_i,u_j\rangle$.
As $u_j$ is in at most two planes,
there is a hyperplane $H$ of $\oR^4$ such that $H\cap P$ is equal to the line
$\langle u_j\rangle$,
and such that all neighbours $t$ of $j$ satisfy $u_t\in H\cup P$.
As $j$ has a neighbour $i$ with $u_i$ at one side of $H$, it also has a
neighbour $k$ with $u_k$ at the other side of $H$, by \rf{27no15d}(ii).
So $u_k\in P$.
Repeating this for $jk$ instead of $ij$, and iterating,
we obtain an infinite walk in $G$, and hence a circuit $C$.
This circuit satisfies the conditions of
Lemma \ref{18no15i}, giving the matrix $A_P$.
\bx

\prop{26no15b}{
Let $G,M,u$ satisfy \rf{26no15a}, with $G$ a $4$-connected flat graph and $d=4$.
Then there exist distinct planes $P_1,P_2,P_3,P_4\subseteq|G|$ with
$P_1\cap P_2\cap P_3\neq\{0\}$ and $P_1\cap P_4=\{0\}$.
}

\pf
By Proposition \ref{26no15bx}, there exist planes $P,R\subseteq|G|$ with
$P\cap R=\{0\}$.
If planes as required do not exist, we can apply Lemma \ref{18no15h} both to $P$
and to $R$.
Consider the graphs $G_P$ and $G_R$ and the
matrices $A_P$ and $A_R$ as in Lemma \ref{18no15h}.
From these we can construct a graph $G'=G_P\cup G_R$ and
a matrix $A:=A_P+A_R$ (where we may assume that $A_P$ and $A_R$ are 0 outside $P$ and $R$
respectively) satisfying the
conditions of Lemma \ref{18no15j}, however with $\lambda^-(A)\geq 2$
(as $P\cap R=\{0\}$), contradicting Lemma \ref{18no15j}.
\bx

\sectz{Theorem and proof}

Having all ingredients, the proof of the theorem now is easy.

\thmnn{
Let $G$ be a $4$-connected flat graph.
Then each well-signed $G$-matrix $M$ with one negative eigenvalue has the
Strong Arnold Property.
}

\pf
Suppose $M$ has not the Strong Arnold Property.
Let $d:=\corank(M)$ and let $u:V(G)\to\oR^d$ be the nullspace embedding defined by $M$.
By [6], $d\leq 4$.
Then Propositions \ref{25no15a}, \ref{18no15a}, and \ref{26no15b} imply $d\leq 3$.

Let $Q$ be a homogeneous quadric in $\oR^d$ with $|G|\subseteq Q$.
By Proposition \ref{18no15a}, $Q$ is not the union of two hyperplanes.
This implies that $d=3$ and that $Q$, and hence $|G|$, contains no plane.
So if $i$ and $j$ are adjacent, then $\dim\langle u_i,u_j\rangle\leq 1$.
Let $H$ be a spanned plane.
Then by \rf{27no15d}(i), $|G|\setminus H$ has at most two components.
Hence it is contained in the union of at most two lines.
So $|G|$ is contained in the union of two planes, a contradiction.
\bx

\bigskip
\noindent
{\em Acknowledgements.}
The research leading to these results has received funding from the European Research Council
under the European Union's Seventh Framework Programme (FP7/2007-2013) / ERC grant agreement
n$\mbox{}^{\circ}$ 339109.

\section*{References}\label{REF}
{\small
\begin{itemize}{}{
\setlength{\labelwidth}{4mm}
\setlength{\parsep}{0mm}
\setlength{\itemsep}{0mm}
\setlength{\leftmargin}{5mm}
\setlength{\labelsep}{1mm}
}
\item[\mbox{\rm[1]}] Y. Colin de Verdi\`ere, 
Sur un nouvel invariant des graphes et un crit\`ere de planarit\'e,
{\em Journal of Combinatorial Theory, Series B} 50 (1990) 11--21
[English translation:
On a new graph invariant and a criterion for planarity,
in: {\em Graph Structure Theory} (N. Robertson, P. Seymour, eds.),
American Mathematical Society,
Providence, Rhode Island, 1993, pp. 137--147].

\item[\mbox{\rm[2]}] H. van der Holst, 
A short proof of the planarity characterization of Colin de Verdi\`ere,
{\em Journal of Combinatorial Theory, Series B} 65 (1995) 269--272.

\item[\mbox{\rm[3]}] H. van der Holst, 
{\em Topological and Spectral Graph Characterizations},
Ph.D.\ Thesis, Universiteit van Amsterdam,
Amsterdam, 1996.

\item[\mbox{\rm[4]}] H. van der Holst, L. Lov\'asz, A. Schrijver, 
The Colin de Verdi\`ere graph parameter,
in: {\em Graph Theory and Combinatorial Biology}
(L. Lov\'asz, et al., eds),
J\'anos Bolyai Mathematical Society, Budapest, 1999, pp. 29--85.

\item[\mbox{\rm[5]}] L. Lov\'asz, 
Steinitz representations and the Colin de Verdi\`ere number,
{\em Journal of Combinatorial Theory, Series B} 82 (2001) 223--236.

\item[\mbox{\rm[6]}] L. Lov\'asz, A. Schrijver, 
A Borsuk theorem for antipodal links and a spectral characterization
of linklessly embeddable graphs,
{\em Proceedings of the American Mathematical Society} 126 (1998) 1275--1285.

\item[\mbox{\rm[7]}] L. Lov\'asz, A. Schrijver, 
On the null space of a Colin de Verdi\`ere matrix,
{\em Annales de l'Institut Fourier}
(Universit\'e de Grenoble) 49 (1999) 1017--1026.

\item[\mbox{\rm[8]}] L. Lov\'asz, A. Schrijver, 
Nullspace embeddings for outerplanar graphs,
pre\-print, 2015.

\item[\mbox{\rm[9]}] N. Robertson, P.D. Seymour, 
Graph minors. {XX}. Wagner's conjecture,
{\em Journal of Combinatorial Theory, Series B} 92 (2004) 325--357.

\item[\mbox{\rm[10]}] N. Robertson, P. Seymour, R. Thomas, 
Sachs' linkless embedding conjecture,
{\em Journal of Combinatorial Theory, Series B} 64 (1995) 185--227.

\end{itemize}
}

\end{document}